\documentclass{article}
\usepackage{amssymb}
\usepackage{amsmath}

\newtheorem{theorem}{Theorem}

\newtheorem{corollary}[theorem]{Corollary}

\newtheorem{lemma}[theorem]{Lemma}

\newtheorem{proposition}[theorem]{Proposition}

\newtheorem{theorem 1}[theorem]{Theorem}
\newtheorem{theorem 2}[theorem]{Theorem}
\newtheorem{theorem 3}[theorem]{Theorem}
\newtheorem{theorem 5}[theorem]{Theorem}
\newtheorem{theorem 6}[theorem]{Theorem}
\newtheorem{theorem 7}[theorem]{Theorem}
\newenvironment{proof6}[1][Proof of theorem 6]{\noindent\textbf{#1.} }{\ \rule{0.5em}{0.5em}}
\newenvironment{proof7}[1][Proof of theorem 7]{\noindent\textbf{#1.} }{\ \rule{0.5em}{0.5em}}
\newenvironment{proofp}[1][Proof of the proposition 4]{\noindent\textbf{#1.} }{\ \rule{0.5em}{0.5em}}
\newenvironment{proofl}[1][Proof of the Lemma 8]{\noindent\textbf{#1.} }{\ \rule{0.5em}{0.5em}}
\newenvironment{proofc}[1][Proof of the corollary 5]{\noindent\textbf{#1.} }{\ \rule{0.5em}{0.5em}}
\begin{document}
	\title{A generalization of the boundedness of certain integral operators in variable Lebesgue spaces}
\author{Marta Urciuolo - Lucas Vallejos}
\maketitle

\begin{abstract}
	Let $A_{1},...A_{m}$ be a $n\times n$ invertible matrices. Let $0 \leq \alpha<n$ and $0<\alpha_{i}<n$ such that $\alpha_1 + ... + \alpha_m = n- \alpha$. We define%
	\begin{equation*}
	T_{\alpha}f(x)=\int \frac{1}{\left\vert x-A_{1}y\right\vert ^{\alpha
			_{1}}...\left\vert x-A_{m}y\right\vert ^{\alpha _{m}}}f(y)dy.
	\end{equation*}%
	In \cite{U-V} we obtained the boundedness of this operator from $L^{p(.)}(%
	\mathbb{R}^{n})$ into $L^{q(.)}(\mathbb{R}^{n})$ for $\frac{1}{q(.)}=\frac{1%
	}{p(.)}-\frac{\alpha }{n},$ in the case that $A_{i}$ is a power of certain
	fixed matrix $A~\ $and for exponent functions $p$ satisfying log-Holder
	conditions and $p(Ay)=p(y),$ $y\in \mathbb{R}^{n}$ $.$ We will show now that
	the hypothesis on $p$, in certain cases, is necessary for the boundedness of $T_{\alpha}$ and we
	also prove the result for more general matrices $A_{i}.$
\footnote{Partially supported by CONICET and SECYTUNC}
\footnote{Math. subject classification: 42B25, 42B35.}
\footnote{Key words: Variable Exponents, Fractional Integrals.}
\end{abstract}

\section{Introduction\protect\bigskip}

Given a measurable function $p(.):\mathbb{R}^{n}\rightarrow \left[ 1,\infty
\right) ,$ let $L^{p(.)}(\mathbb{R}^{n})$ be the Banach space of measurable
functions $f$ on $\mathbb{R}^{n}$ such that for some $\lambda >0,$ 
\begin{equation*}
\int \left( \frac{\left\vert f(x)\right\vert }{\lambda }\right)
^{p(x)}dx<\infty ,
\end{equation*}%
with norm 
\begin{equation*}
\left\Vert f\right\Vert _{p(.)}=\inf \left\{ \lambda >0:\int \left( \frac{%
\left\vert f(x)\right\vert }{\lambda }\right) ^{p(x)}dx\leq 1\right\} .
\end{equation*}%
These spaces are known as \textit{variable exponent spaces} and are a
generalization of the classical Lebesgue spaces $L^{p}(\mathbb{R}^{n}).$
They have been widely studied lately. See for example \cite{CCF}, 
\cite{C-F-N} and \cite{DHHR}. The first step was to determine sufficient
conditions on $p(.)$ for the boundedness on $L^{p(.)}$ of the Hardy
Littlewood maximal operator 
\begin{equation*}
\mathcal{M}f(x)=\sup\limits_{B}\frac{1}{\left\vert B\right\vert }%
\int_{B}\left\vert f(y)\right\vert dy,
\end{equation*}%
where the supremun is taken over all balls $B$ containing $x.$ Let $%
p_{-}=ess\inf $\textit{\ }$p(x)$ and let $p_{+}=ess\sup $ $p(x).$ In \cite
{C-F-N} , D. Cruz Uribe, A. Fiorenza and C. J. Neugebauer proved the following result.

%\textbf{Theorem} \textit{
\begin{theorem 1}
Let $p(.):\mathbb{R}^{n}\rightarrow \left[
1,\infty \right) $ \textit{\ be such that }$1<p_{-}\leq p_{+}<\infty .$ 
\textit{Suppose further that }$p(.)$\textit{\ satisfies } 
\begin{equation}
\left\vert p(x)-p(y)\right\vert \leq \frac{c}{-\log \left\vert
x-y\right\vert },\text{ }\left\vert x-y\right\vert <\frac{1}{2},
\label{log local}
\end{equation}%
\textit{and } 
\begin{equation}
\left\vert p(x)-p(y)\right\vert \leq \frac{c}{\log \left( e+\left\vert
x\right\vert \right) },\text{ }\left\vert y\right\vert \geq \left\vert
x\right\vert .  \label{log fuera}
\end{equation}%
\textit{Then the Hardy Littlewood maximal operator is bounded on }$L^{p(.)}(%
\mathbb{R}^{n}).$
\end{theorem 1}

We recall that a weight $\omega$ is a locally integrable and non negative
function. The Muckenhoupt class $\mathcal{A}_p$, $1<p<\infty$, is defined as the class
of weights $\omega$ such that

\begin{equation*}
\sup_{Q} \left[ \left( \frac{1}{\lvert Q \rvert} \displaystyle \int_Q \omega
\right) \left( \frac{1}{\lvert Q \rvert} \displaystyle \int_Q \omega^{-\frac{%
1}{p-1}} \right)^{p-1} \right] < \infty,
\end{equation*}
where $Q$ is a cube in $\mathbb{R}^n$.

For $p=1$, $\mathcal{A}_1$ is the class of weights $\omega$ satisfying that there
exists $c>0$ such that

\begin{equation*}
\mathcal{M} \omega (x) \leq c \omega (x) \ a.e. \ x \in \mathbb{R}^n.
\end{equation*}

We denote $\left[ \omega \right]_{\mathcal{A}_1}$ the infimum of the constant $c$ such
that $\omega$ satisfies the above inequation.

In \cite{M-W}, B. Muckenhoupt y R.L. Wheeden define $\mathcal{A}(p,q),$ $1<p<\infty$ and $%
1<q<\infty$, as the class of weights $\omega$ such that

\begin{equation*}
\sup_{Q} \displaystyle\left[ \left( \dfrac{1}{\lvert Q \rvert} \displaystyle%
\int_{Q}\omega(x)^q dx \right)^{\frac{1}{q}}\left( \dfrac{1}{\lvert Q \rvert}
\displaystyle\int_{Q} \omega(x)^{-p^{\prime }}dx \right)^{\frac{1}{p^{\prime
}}}\displaystyle\right] \ < \ \infty.
\end{equation*}

When $p=1$, $\omega \in \mathcal{A}(1,q)$ if only if

\begin{equation*}
\sup_{Q} \displaystyle\left[ \lVert \omega^{-1} \chi_Q \rVert_{\infty}\left( 
\dfrac{1}{\lvert Q \rvert} \displaystyle\int_Q \omega(x)^{q}dx \right)^{%
\frac{1}{q}}\displaystyle\right] \ < \ \infty.
\end{equation*}
Let $0\leq \alpha <n$. For $1\leq i\leq m$, let $0<\alpha _{i}<n,$ be such
that 
\begin{equation*}
\alpha _{1}+...+\alpha _{m}=n-\alpha .
\end{equation*}%
Let $T_{\alpha}$ be the integral operator given by

\begin{equation}
T_{\alpha}f\left( x\right) =\int k\left( x,y\right) f\left( y\right) dy
\label{TA}
\end{equation}%
where 
\begin{equation*}
k\left( x,y\right) =\frac{1}{\left\vert x-A_{1}y\right\vert ^{\alpha _{1}}}%
...\frac{1}{\left\vert x-A_{m}y\right\vert ^{\alpha _{m}}},  \label{K}
\end{equation*}%
and where the matrices $A_{i}$ are certain invertible matrices such that $%
A_{i}-A_{j}$ is invertible for $i\neq j,$ $1\leq i,j\leq m.$

In the paper \cite{Ri-U} the authors studied this kind of integral operators
and they obtained weighted $\left( p,q\right) $ estimates, $\frac{1}{q}=%
\frac{1}{p}-\frac{\alpha }{n},$ for weights $w\in A(p,q)$ such that $%
w(A_{i}x)\leq cw(x).$ In \cite{U-V} we use extrapolation techniques to
obtain $p(.)-q(.)$ and weak type estimates, in the case where $A_{i}=A^{i},$
and $A^{N}=I,$ for some $N\in \mathbb{N}$. This technique allows us to
replace the log-H\"{o}lder conditions about the exponent $p(\cdot )$ by a
more general hypothesis concerning the boundeness of the maximal function $%
\mathcal{M}$. We obtain the following results

\begin{theorem 2}
Let $A$ be an invertible matrix such that $%
A^{N}=I,$ for some $N\in \mathbb{N}$, let $T_{\alpha}$ be the integral operator given by ($\ref{TA}$), where $A_{i}=A^{i}$ and such that $A_{i}-A_{j}$ is invertible
for $i\neq j$, $1\leq i,j\leq m$. Let $p:\mathbb{R}^{n}\longrightarrow \left[
1,\infty \right) $ be such that $1<p_{-}\leq p_{+}<\frac{n}{\alpha }$ and
such that $p(Ax)=p(x)$ a.e. $x\in \mathbb{R}^{n}$. Let $q(\cdot )$ be
defined by $\frac{1}{p(x)}-\frac{1}{q(x)}=\frac{\alpha }{n}.$ If the maximal
operator $\mathcal{M}$ is bounded on $L^{\left( \frac{n-\alpha p_{-}}{np_{-}}%
q(.)\right) ^{\prime }}$ then $T$ is bounded from $L^{p(.)}\left( \mathbb{R}%
^{n}\right) $ into $L^{q(.)}(\mathbb{R}^{n}).$
\end{theorem 2}

%\textbf{Theorem} \textbf{B}
\begin{theorem 3}
Let $A$ be an invertible matrix such that $%
A^{N}=I,$ for some $N\in \mathbb{N}$, let $T_{\alpha}$ be the integral operator given
by ($\ref{TA}$), where $A_{i}=A^{i}$ and such that $A_{i}-A_{j}$ is invertible
for $i\neq j$, $1\leq i,j\leq m$. Let $p:\mathbb{R}^{n}\longrightarrow \left[
1,\infty \right) $ be such that $1\leq p_{-}\leq p_{+}<\frac{n}{\alpha }$
and such that $p(Ax)=p(x)$ a.e. $x\in \mathbb{R}^{n}$. Let $q(\cdot )$ be
defined by $\frac{1}{p(x)}-\frac{1}{q(x)}=\frac{\alpha }{n}.$ If the maximal
operator $\mathcal{M}$ is bounded on $L^{\left( \frac{n-\alpha p_{-}}{np_{-}}%
q(.)\right) ^{\prime }}$ then there exists $c>0$ such that 
\begin{equation*}
\left\Vert t\chi _{\left\{ x:T_{\alpha }f(x)>t\right\} }\right\Vert
_{q(.)}\leq c\left\Vert f\right\Vert _{p(.)}.
\end{equation*}
\end{theorem 3}

We also showed that this technique applies in the case when each of the
matrices $A_{i}$ is either a power of an orthogonal matrix $A$ or a power of 
$A^{-1}$.

In this paper we will prove that these theorems generalize to any invertible
matrices $A_{1},...,A_{m}$ such that $A_{i}-A_{j}$ is invertible for $i\neq
j,$ $1\leq i,j\leq m.$ We will also show, in some cases, that the
condition $p(A_{i}x)=p(x),$ $x\in \mathbb{R}^{n}$ is necessary to obtain $%
p(.)-q(.)$ boundedness.

\section{Necessary conditions on $p$}

Let $A$ be a $n\times n$ invertible matrix and let $%
0<\alpha <n$. We define%
\begin{equation*}
T_{A}f(x)=\int \frac{1}{\left\vert x-Ay\right\vert ^{n-\alpha }}f(y)dy.
\end{equation*}

\begin{proposition}
Let $A$ be a $n\times n$ invertible matrix. Let $p:\mathbb{R}^{n}\rightarrow %
\left[ 1,\infty \right) $ be a measurable function such that $p$ is
continuos at $y_{0}$ and at $Ay_{0}$ for some $y_{0}\in \mathbb{R}^{n}.$ If $%
p(Ay_{0})> p(y_{0})$ then there exists $f\in L^{p(.)}(\mathbb{R}^{n})$
such that $T_{A}f\notin L^{q(.)}(\mathbb{R}^{n})$ for $\frac{1}{q(.)}=\frac{1%
}{p(.)}-\frac{\alpha }{n}.$
\end{proposition}

\begin{proofp}
Since $p$ is continuos at $y_{0},$ there exists \ ball \thinspace $%
B=B(y_{0},r)$ such that $p(y)\thicksim p(y_{0})$ for $x\in B.$ We suppose $%
p(y_{0})<p(Ay_{0}).$ In this case we take 
\begin{equation*}
f(y)=\frac{\chi _{B}(y)}{\left\vert y-y_{0}\right\vert ^{\beta }},
\end{equation*}%
for certain $\beta <\frac{n}{p(y_{0})}$ that will be chosen later. We will
show that, for certain $\beta ,$ $f\in L^{p(.)}(\mathbb{R}^{n})$ but $%
T_{A}f\notin L^{q(.)}(\mathbb{R}^{n}).$ Indeed, 
\begin{equation*}
T_{A}f(x)=\int \frac{1}{\left\vert x-Ay\right\vert ^{n-\alpha }}%
f(y)dy=\int_{B}\frac{1}{\left\vert x-Ay\right\vert ^{n-\alpha }\left\vert
y-y_{0}\right\vert ^{\beta }}dy,
\end{equation*}%
so 
\begin{equation*}
\int \left( T_{A}f(x)\right) ^{q(x)}dx=\int \left( \int_{B}\frac{1}{%
\left\vert x-Ay\right\vert ^{n-\alpha }\left\vert y-y_{0}\right\vert ^{\beta
}}dy\right) ^{q(x)}dx
\end{equation*}%
\begin{equation*}
\geq \int_{B\left( Ay_{0},\varepsilon \right) }\left( \int_{B}\frac{1}{%
\left\vert x-Ay\right\vert ^{n-\alpha }\left\vert y-y_{0}\right\vert ^{\beta
}}dy\right) ^{q(x)}dx
\end{equation*}%
\begin{equation*}
\geq \int_{B\left( Ay_{0},\varepsilon \right) }\left( \int_{B\cap \left\{
y:\left\vert Ay-Ay_{0}\right\vert <\left\vert Ay_{0}-x\right\vert \right\} }%
\frac{1}{\left\vert x-Ay\right\vert ^{n-\alpha }\left\vert
y-y_{0}\right\vert ^{\beta }}dy\right) ^{q(x)}dx
\end{equation*}%
Now, we denote by $M=\left\Vert A\right\Vert =\sup\limits_{\left\Vert
y\right\Vert =1}\left\vert Ay\right\vert $. Now for $\varepsilon <Mr$ and $%
x\in B\left( Ay_{0},\varepsilon \right) $, $B(y_{0},\frac{1}{M}\left\vert
Ay_{0}-x\right\vert )\subset B\cap \left\{ y:\left\vert Ay-Ay_{0}\right\vert
<\left\vert Ay_{0}-x\right\vert \right\} .$ Indeed, $\left\vert
y-y_{0}\right\vert \leq \frac{1}{M}\left\vert Ay_{0}-x\right\vert \leq \frac{%
1}{M}\varepsilon \leq r$ and $\left\vert Ay-Ay_{0}\right\vert \leq
M\left\vert y-y_{0}\right\vert \leq \left\vert Ay_{0}-x\right\vert ,$ so 
\begin{equation*}
\geq \int_{B\left( Ay_{0},\varepsilon \right) }\left( \int_{B(y_{0},\frac{1}{%
M}\left\vert Ay_{0}-x\right\vert )}\frac{1}{\left\vert x-Ay\right\vert
^{n-\alpha }\left\vert y-y_{0}\right\vert ^{\beta }}dy\right) ^{q(x)}dx,
\end{equation*}%
also, for $y\in B(y_{0},\frac{1}{M}\left\vert Ay_{0}-x\right\vert )$%
\begin{equation*}
\left\vert x-Ay\right\vert \leq \left\vert x-Ay_{0}\right\vert +\left\vert
Ay_{0}-Ay\right\vert \leq \left\vert x-Ay_{0}\right\vert +M\left\vert
y_{0}-y\right\vert \leq 2\left\vert x-Ay_{0}\right\vert ,
\end{equation*}%
so%
\begin{equation*}
\geq \int_{B\left( Ay_{0},\varepsilon \right) }\left( \frac{1}{2^{n-\alpha
}\left\vert x-Ay_{0}\right\vert ^{n-\alpha }}\right) ^{q(x)}\left(
\int_{B(y_{0},\frac{1}{M}\left\vert Ay_{0}-x\right\vert )}\frac{1}{%
\left\vert y-y_{0}\right\vert ^{\beta }}dy\right) ^{q(x)}dx
\end{equation*}%
\begin{equation*}
=\int_{B\left( Ay_{0},\varepsilon \right) }\left( \frac{1}{2^{n-\alpha
}\left\vert x-Ay_{0}\right\vert ^{n-\alpha }}\right) ^{q(x)}\left(
c\left\vert Ay_{0}-x\right\vert ^{-\beta +n}\right) ^{q(x)}dx
\end{equation*}%
\begin{equation*}
=\int_{B\left( Ay_{0},\varepsilon \right) }\left( \frac{c}{2^{n-\alpha
}\left\vert x-Ay_{0}\right\vert ^{\beta -\alpha }}\right) ^{q(x)}dx.
\end{equation*}%
Now, since $q(Ay_{0})>q(y_{0}),$ $q(Ay_{0})-\gamma >q(y_{0})$ for $\gamma =%
\frac{q(Ay_{0})-q(y_{0})}{2}$ We observe that if $\frac{1}{q(y_{0})}=\frac{1%
}{p(y_{0})}-\frac{\alpha }{n},$ for $\beta _{0}=\frac{n}{p(y_{0})},$ $\left(
\beta _{0}-\alpha \right) q(y_{0})=\left( \frac{n}{p(y_{0})}-\alpha \right)
q(y_{0})=n,$ so since $q(Ay_{0})-\gamma >q(y_{0}),$ we obtain that $\left( 
\frac{n}{p(y_{0})}-\alpha \right) \left( q(Ay_{0})-\gamma \right) >n$. and
still $\left( \beta -\alpha \right) \left( q(Ay_{0})-\gamma \right) >n$ for $%
\beta =\frac{n}{p(y_{0})}-\frac{1}{2}\left( \frac{n}{p(y_{0})}-\left( \alpha
+\frac{n}{q(Ay_{0})-\gamma )}\right) \right) $. So $\beta =\frac{n}{p(y_{0})}%
(1-\delta )$ for some $\delta >0.$ Since $q$ is continuos, we chose $%
\varepsilon $ so that, for $x\in B\left( Ay_{0},\varepsilon \right) ,$ $%
q(x)>q(Ay_{0})-\gamma .$ and $\frac{c}{2^{n-\alpha }\left\vert
x-Ay_{0}\right\vert ^{\beta -\alpha }}>1$ so this last integral is bounded
from below by 
\begin{equation*}
c\int_{B\left( Ay_{0},\varepsilon \right) }\left( \frac{1}{\left\vert
x-Ay_{0}\right\vert ^{\beta -\alpha }}\right) ^{q(Ay_{0})-\gamma }dx=\infty .
\end{equation*}%
For this $\beta $ we chose $r$ to obtain that the ball $B=B(y_{0},r)\subset
\left\{ y:p(y)<\frac{p(y_{0})}{1-\delta }\right\} $. In this way we obtain
that $f\in L^{p(.)}(\mathbb{R}^{n})$ but $T_{A}f\notin L^{q(.)}(\mathbb{R}%
^{n}).$ \newline\newline
\end{proofp}

\begin{corollary}
	If $A^N=I$ for some $N \in \mathbb{N}$, $p$ is continuos and $T_A$ bounded is from $L^{p(.)}$ into $L^{q(.)}$, then $p(Ay)=p(y)$ for all $y \in \mathbb{R}^n$
\end{corollary}

\begin{proofc}
We suposse that $p(Ay_0) < p(y_0)$. Since $p$ is continuos in $y_0$, by the last proposition,
\begin{equation*}
p(Ay_0) < p( y_0) = p(A^N y_0) \leq p(A^{N-1}y_0) \leq ... \leq p(Ay_0) = p(Ay_0) 
\end{equation*}
which is a contradiction.
\end{proofc}

\section{The main results}

Given $0\leq \alpha <n,$ we recall that we are studying fractional type
integral operators of the form 
\begin{equation}
T_{\alpha}f\left( x\right) =\int k\left( x,y\right) f\left( y\right) dy,
\end{equation}%
with a kernel 
\begin{equation*}
k(x,y)=\frac{1}{\left\vert x-A_{1}y\right\vert ^{\alpha _{1}}}...\frac{1}{%
\left\vert x-A_{m}y\right\vert ^{\alpha _{m}}},
\end{equation*}%
$\alpha _{1}+...+\alpha _{m}=n-\alpha ,$ $0<\alpha _{i}<n.$

\begin{theorem 6}
Let $m\in \mathbb{N},$ let $A_{1},...A_{m}$ be invertibles matrices such that 
$A_{i}-A_{j}$ is invertible for $i\neq j$, $1\leq i,j\leq m$. Let $T_{\alpha}$ be the
integral operator given by (\ref{TA}), let $p:\mathbb{R}^{n}\longrightarrow %
\left[ 1,\infty \right) $ be such that $1<p_{-}\leq p_{+}<\frac{n}{\alpha }$
and such that $p(A_{i}x)=p(x)$ a.e. $x\in \mathbb{R}^{n},$ $1\leq i\leq m$.
Let $q(\cdot )$ be defined by $\frac{1}{p(x)}-\frac{1}{q(x)}=\frac{\alpha }{n%
}.$ If the maximal operator $\mathcal{M}$ is bounded on $L^{\left( \frac{%
n-\alpha p_{-}}{np_{-}}q(.)\right) ^{\prime }}$ then $T_{\alpha}$ is bounded from $%
L^{p(.)}\left( \mathbb{R}^{n}\right) $ into $L^{q(.)}(\mathbb{R}^{n}).$
\end{theorem 6}

\begin{theorem 7}
Let $m\in \mathbb{N},$ let $A_{1},...A_{m}$ be invertibles matrices such that 
$A_{i}-A_{j}$ is invertible for $i\neq j$, $1\leq i,j\leq m$. Let $T_{\alpha}$ be the
integral operator given by (\ref{TA}), let $p:\mathbb{R}^{n}\longrightarrow %
\left[ 1,\infty \right) $ be such that $1\leq p_{-}\leq p_{+}<\frac{n}{%
\alpha }$ and such that $p(A_{i}x)=p(x)$ a.e. $x\in \mathbb{R}^{n},$ $1\leq
i\leq m$. Let $q(\cdot )$ be defined by $\frac{1}{p(x)}-\frac{1}{q(x)}=\frac{%
\alpha }{n}.$ If the maximal operator $\mathcal{M}$ is bounded on $L^{\left( 
\frac{n-\alpha p_{-}}{np_{-}}q(.)\right) ^{\prime }}$ then there exists $c>0$
such that 
\begin{equation*}
\left\Vert t\chi _{\left\{ x:T_{\alpha }f(x)>t\right\} }\right\Vert
_{q(.)}\leq c\left\Vert f\right\Vert _{p(.)},
\label{td}
\end{equation*}
\end{theorem 7}

\begin{lemma}
	If $f \in L^{1}_{loc}(\mathbb{R}^n)$ and $A$ an invertible $n \times n$ matrix then
	\begin{equation*}
	\mathcal{M} (f \circ A)(x) \leq c (\mathcal{M}(f) \circ A)(x).
	\end{equation*}
\end{lemma}

\begin{proofl}
	Indeed, $\mathcal{M}(f \circ A)=\sup\limits_{B} \dfrac{1}{\lvert B \rvert} \displaystyle\int_{B}\left\vert (f\circ A)(y)\right\vert dy$, where the suppremun is taken over all balls $B$ containing $x$. By a change of variable we see that,\newline 
	\begin{equation*}
	\dfrac{1}{\lvert B \rvert}\displaystyle\int_{B}\left\vert (f\circ A)(y)\right\vert dy = \lvert det(A^{-1}) \rvert \dfrac{1}{\lvert B \rvert}\displaystyle\int_{A(B)}\left\vert f(z)\right\vert dz,  
	\end{equation*}
	where $A(B)=\{Ay:y \in B\}$. Now, if $y \in B=B(x_0 , r)$ then $\lvert Ay-Ax_0 \rvert \leq M \lvert y-x_0 \rvert \leq Mr$, where $M=\lVert A \rVert$. That is $Ay \in \widetilde{B} = B(Ax_0 , Mr)$, So\newline
	\begin{equation*}
	\leq \dfrac{M^{n} \lvert det(A^{-1}) \rvert}{\lvert \widetilde{B} \rvert} \displaystyle \int_{\widetilde{B}} f(z) dz
	\end{equation*}
	\begin{equation*}
	\leq M^{n} \lvert det(A^{-1}) \rvert \mathcal{M}f(Ax).
	\end{equation*}
	Therefore we obtain that, 
	\begin{equation*}
	\mathcal{M} (f \circ A) \leq c (\mathcal{M}(f) \circ A),
	\end{equation*} 
	with $c=M^{n} \lvert det(A^{-1}) \rvert$.
	%It will be enough to see that $\mathcal{R}h^{\frac{1}{q_0}} \circ A_j \in A_1$ for all $j=1...m$. 
\end{proofl}

\section{Proofs of the main results}
\begin{proof6}
In the paper \cite{Ri-U} the authors obtain an estimate of the form
\begin{equation}
\int \left( T_{\alpha }f\right) ^{p}(x)w(x)dx\leq
c \sum\limits_{j=1}^{m}\int \left( \mathcal{M}_{\alpha }f\right)
^{p}(x)w(A_{j}x)dx,  \label{T-M}
\end{equation}%
for any $w\in \mathcal{A}_{\infty }$ and $0<p<\infty $ (See the last lines of page 454
in \cite{Ri-U})$.$ We denote $\widetilde{q}(.)=\frac{q(.)}{q_{0}},$ we
define an iteration algorithm on $L^{\widetilde{q}(.)^{\prime }}$ by 
\begin{equation}
\mathcal{R}h(x)=\sum_{k=0}^{\infty }\frac{\mathcal{M}%
	^{k}h(x)}{2^{k}\left\Vert \mathcal{M}\right\Vert _{\widetilde{q}%
		(.)^{\prime }{}}^{k}},
\label{Rh}
\end{equation}%
where, for $k\geq 1$, $\mathcal{M}^{k}$ denotes $k$ iteration of the
maximal operator $\mathcal{M}$ and $\mathcal{M}^{0}\left( h\right)
=\left\vert h\right\vert .$ We will check that\newline
$a)$ For all $x\in \mathbb{R}^{n}, \left\vert
h(x)\right\vert \leq \mathcal{R}h(x),$\newline
$b)$ For all $j:1,...,m, \left\Vert \mathcal{R}h\circ A_j \right\Vert _{\widetilde{%
		q}(.)^{\prime }}\leq c\left\Vert h\right\Vert _{\widetilde{q}(.)^{\prime }},$\newline
$c)$ For all $j:1, ..., m, (\mathcal{R}h \circ A_j)^{\frac{1}{q_0}} \in \mathcal{A}(p_{-},q_0)$.\newline\newline
Indeed, $a)$ is evident from the definition $b)$ is verified by the
following,
\begin{equation*}
\lVert \mathcal{R}h \circ A_j \rVert_{\widetilde{q}(.)^{'}} \leq \displaystyle \sum \limits_{k=0}^{\infty} \dfrac{\lVert \mathcal{M}^k h \circ A_j\rVert_{\widetilde{q}(.)^{'}}}{2^k \lVert \mathcal{M} \rVert^{k}_{\widetilde{q}(.)^{'}} }
\end{equation*}
and \newline\newline
\begin{equation*}
\lVert \mathcal{M}^k h \circ A_j\rVert_{\widetilde{q}(.)^{'}} = \inf \left\{\lambda > 0 : \displaystyle \int_{\mathbb{R}^{n}} \left(\dfrac{\mathcal{M}^k h (A_j x)}{\lambda}\right)^{\widetilde{q}(x)^{'}} dx \leq 1 \right\}
\end{equation*}
But, by a change of variable and using the hypothesis on the exponent,\newline\newline
\begin{equation*}
\displaystyle \int_{\mathbb{R}^{n}} \left(\dfrac{\mathcal{M}^k h (A_j x)}{\lambda}\right)^{\widetilde{q}(x)^{'}} dx = \lvert det(A_{j}^{-1}) \rvert \displaystyle \int_{\mathbb{R}^{n}} \left(\dfrac{\mathcal{M}^k h (y)}{\lambda}\right)^{\widetilde{q}(A_{j}^{-1}y)^{'}} dy,
\end{equation*}\newline
put D=$max \left\{ \lvert det(A_{j}^{-1}) \rvert, j=1...m \right\}$,\newline
\begin{equation}
\leq D \displaystyle \int_{\mathbb{R}^n} \left( \dfrac{\mathcal{M}^k h (y)}{\lambda} \right)^{\widetilde{q}^{'}(y)} dy
\end{equation} 
If $D \leq 1$,
%That is, 
\begin{equation*}
\lVert \mathcal{M}^k h \circ A_j\rVert_{\widetilde{q}(.)^{'}} \leq  \lVert \mathcal{M}^k h\rVert_{\widetilde{q}(.)^{'}}
\end{equation*}
So,
\begin{equation*}
	\lVert \mathcal{R}h \circ A_j \rVert_{\widetilde{q}(.)^{'}} \leq \sum_{k=0}^{\infty }\frac{\lVert \mathcal{M}%
		^{k}h(x) \rVert_{\widetilde{q}(.)^{'}}}{2^{k}\left\Vert \mathcal{M}\right\Vert _{\widetilde{q}%
			(.)^{\prime }{}}^{k}} \leq \lVert h \rVert_{\widetilde{q}(.)^{'}} \sum_{k=0}^{\infty} \dfrac{1}{2^k}=2 \lVert h \rVert_{\widetilde{q}(.)^{'}}
\end{equation*}

If $D > 1$ then from (7) it is follows that
\begin{equation*}
D \displaystyle \int_{\mathbb{R}^n} \left( \dfrac{\mathcal{M}^k h (y)}{\lambda} \right)^{\widetilde{q}^{'}(y)} dy = \displaystyle \int_{\mathbb{R}^n} \left( \dfrac{M^k h(y)}{\lambda C^{\frac{1}{\widetilde{q}(.)^{'}}}} \right)^{\widetilde{q}(.)^{'}} dy
\end{equation*} 
and $D=\frac{1}{C}$ where $C=min\{\lvert det(A_j)\rvert,  j=1...m \}$. So,
\begin{equation*}
\leq \displaystyle \int_{\mathbb{R}^n} \left( \dfrac{M^k h(y)}{\lambda C^{\frac{1}{(\widetilde{q}^{'})_{-}}}} \right)^{\widetilde{q}(.)^{'}} dy
\end{equation*}
That is,
\begin{equation*}
\displaystyle \int_{\mathbb{R}^{n}} \left(\dfrac{\mathcal{M}^k h (A_j x)}{\lambda}\right)^{\widetilde{q}(x)^{'}} dx \leq \displaystyle \int_{\mathbb{R}^n} \left( \dfrac{M^k h(x)}{\lambda C^{\frac{1}{(\widetilde{q}^{'})_{-}}}} \right)^{\widetilde{q}(.)^{'}} dx.
\end{equation*}
From this last inequality it follows that
\begin{equation*}
\lVert \mathcal{M}^k h \circ A_j\rVert_{\widetilde{q}(.)^{'}} \leq D^{\frac{1}{(\widetilde{q}^{'})_{-}}} \lVert \mathcal{M}^k h\rVert_{\widetilde{q}(.)^{'}}
\end{equation*}
and so $b)$ is verified with $c=2D^{\frac{1}{(\widetilde{q}^{'})_{-}}}$\newline
To see $c)$, by Lemma 8,
\begin{equation*}
\mathcal{M}(\mathcal{R}h^{\frac{1}{q_0}} \circ A_j)(x) \leq c \mathcal{M}(\mathcal{R}h^{\frac{1}{q_0}})(A_j x)
\end{equation*}  
$\mathcal{R}h \in \mathcal{A}_1$ (see \cite{C-F}) implies that $\mathcal{R}h^{\frac{1}{q_0}} \in \mathcal{A}_1$ and so,
\begin{equation*}
\leq c \mathcal{R}h^{\frac{1}{q_0}} (A_j x) = c(\mathcal{R}h^{\frac{1}{q_0}} \circ A_j)(x).
\end{equation*} 
Then $c)$ follows since a weight $\omega \in \mathcal{A}_1$ implies that $\omega \in \mathcal{A}(p_{-},q_0)$. \newline
We now take a bounded function $f$\ with compact support. We will check
later that $\left\Vert T_{\alpha}f\right\Vert _{q(.)}<\infty $ , so as in Theorem
5.24 in \cite{C-F}, 
\begin{equation*}
\left\Vert T_{\alpha}f\right\Vert _{q(.)}^{q_{0}}=\left\Vert \left( T_{\alpha}f\right)
^{q_{0}}\right\Vert _{\widetilde{q}(.)}=C\sup_{\left\Vert h\right\Vert _{%
		\widetilde{q}(.)^{\prime }}=1}\int \left( T_{\alpha}f\right) ^{q_{0}}(x)h(x)dx
\end{equation*}%
\begin{equation*}
\leq C\sup_{\left\Vert h\right\Vert _{\widetilde{q}(.)^{\prime }}=1}\int
\left( T_{\alpha}f\right) ^{q_{0}}(x)\mathcal{R}h(x)dx\leq C\sup_{\left\Vert
	h\right\Vert _{\widetilde{q}(.)^{\prime }}=1}\sum\limits_{j=1}^{m}\int
\left( \mathcal{M}_{\alpha }f\right) ^{q_{0}}(x)\mathcal{R}h(A_{j}x)dx,
\end{equation*}%
\begin{equation*}
\leq C\sup_{\left\Vert h\right\Vert _{\widetilde{q}(.)^{\prime
	}}=1}\sum\limits_{j=1}^{m}\left( \int \left\vert f(x)\right\vert ^{p_{-}}%
\mathcal{R}h^{\frac{p_{-}}{q_{0}}}(A_{j}x)dx\right) ^{\frac{q_{0}}{p_{-}}}
\end{equation*}%
where the last inequality follows since $ \mathcal{R}h^{\frac{1}{q_{0}}} \circ A_i $ are weights in $\mathcal{A}(p_{-},q_{0})$ (by $c)$). We
denote by $\widetilde{p}(.)=\frac{p(.)}{p_{-}}.$ Holder%
%TCIMACRO{\U{b4}}%
%BeginExpansion
\'{}%
%EndExpansion
s inequality, $2)$ and Proposition 2.18 in \cite{C-F} and again the
hypothesis about $A_{i}$ and $p$ give%
\begin{equation*}
\left\Vert T_{\alpha}f\right\Vert _{q(.)}^{q_{0}}\leq C\left\Vert
f^{p_{-}}\right\Vert _{\widetilde{p}(.)}^{\frac{q_{0}}{p_{-}}%
}\sup_{\left\Vert h\right\Vert _{\widetilde{q}(.)^{\prime
	}}=1}\sum\limits_{j=1}^{m}\left\Vert \left( \mathcal{R}h ^{%
	\frac{p_{-}}{q_{0}}} \right) \circ {A_{j}}\right\Vert _{\widetilde{p}(.)^{\prime }}^{\frac{q_{0}}{%
		p_{-}}}
\end{equation*}%
\begin{equation*}
\leq \sup_{\left\Vert h\right\Vert _{\widetilde{q}(.)^{\prime
	}}=1} Cm\left\Vert f\right\Vert _{p(.)}^{q_{0}}\left\Vert
h\right\Vert _{\widetilde{q}(.)^{\prime }} \leq C\left\Vert f\right\Vert _{p(.)}^{q_{0}}.
\end{equation*}%
\newline
Now we show that $\left\Vert T_{\alpha} f\right\Vert _{q(.)}<\infty .$ By Prop. 2.12,
p.19 in \cite{C-F}, it is enough to check that $\rho _{q(.)}\left( T_{\alpha}f\right)
<\infty .$%
\begin{equation*}
\left\vert T_{\alpha}f(x)\right\vert ^{q(x)}\leq \left\vert T_{\alpha}f(x)\right\vert
^{q_{+}}\chi _{\left\{ x:T_{\alpha}f(x)>1\right\} }+\left\vert T_{\alpha}f(x)\right\vert
^{q_{-}}\chi _{\left\{ x:T_{\alpha}f(x)\leq 1\right\} },
\end{equation*}%
now $f$ is bounded and with compact support, so $T_{\alpha}f\in L^{s}(\mathbb{R}^{n})$
for $\frac{n}{n-\alpha }<s<\infty ,$ (see Lemma 2.2 in \cite{Ri-U}) thus $%
\int \left\vert T_{\alpha}f(x)\right\vert ^{q(x)}dx<\infty .$ The theorem follows
since bounded functions with compact support are dense in $L^{p(.)}\left( 
\mathbb{R}^{n}\right) $ (See Corollary 2.73 in \cite{C-F}).	
\end{proof6}\newline

\begin{proof7}
We observe that it is enough to check (\ref{td}) for $f \in L^{\infty}_{c}(\mathbb{R}^n)$.\newline
In \cite{Ri-U} (See page 459) the authors prove that there exists $c>0$ such that,
\begin{equation*}
\sup_{\lambda > 0} \lambda \left( \omega^{q_0} \{ x:\lvert T_{\alpha} f(x) \rvert > \lambda\} \right) ^{\frac{1}{q_0}} \leq
\sup_{\lambda > 0} \lambda \left( \omega^{q_0} \{ x: \sum_{i=1}^{m} \mathcal{M}_{\alpha}f(A^{-1}_i x) > c \lambda \} \right) ^{\frac{1}{q_0}}
\end{equation*}
for all $\omega \in \mathcal{A}_{\infty}$ and $f \in L^{\infty}_{c}(\mathbb{R}^n)$.\newline\newline
Let $F_{\lambda}= \lambda^{q_0} \chi_{ \{ x:\lvert T_{\alpha}f(x) \rvert > \lambda \} }$
the last inequality implies that,
\begin{equation}
\int_{\mathbb{R}^n}F_{\lambda}(x) \omega(x)^{q_0} dx \leq \sup_{\lambda >0} \int_{\mathbb{R}^n}\lambda^{q_0}\chi(x)_{ \{ x: \sum_{i=1}^m \mathcal{M}_{\alpha}f(A^{-1}_i x)> c \lambda \} } \omega(x)^{q_0} dx
\label{DI}
\end{equation}
for some $c>0$ and for all $\omega \in \mathcal{A}_{\infty}$.
Now by proposition 2.18 in \cite{C-F}, if $\widetilde{q}(.)=\frac{q(.)}{q_0}$,
\begin{equation*}
\lVert \lambda \chi_{ \{ x:\lvert T_{\alpha}f(x) \rvert > \lambda \} }  \rVert^{q_0}_{q(.)} = \lVert \lambda^{q_0} \chi_{ \{ x:\lvert T_{\alpha}f(x) \rvert > \lambda \} }  \rVert_{\widetilde{q}(.)}
\end{equation*}
\begin{equation*}
= \lVert F_{\lambda} \rVert_{\widetilde{q}(.)} \leq c \sup_{\left\Vert h\right\Vert_{\widetilde{q}'(.)} = 1} \int_{\mathbb{R}^{n}} F_{\lambda}(x) h(x) dx,
\end{equation*}
Let $\mathcal{R}h$ be define by (\ref{Rh}). We can verify that,\newline

$a)$ $\left\vert h(x)\right\vert \leq \mathcal{R}h(x)$ $x\in \mathbb{R}^{n}$, \newline

$b)$ For all $j:1,...,m, \left\Vert \mathcal{R}h\circ A_j \right\Vert _{\widetilde{%
		q}(.)^{\prime }}\leq c\left\Vert h\right\Vert _{\widetilde{q}(.)^{\prime }}$, \newline

$c)$ For all $j:1, ..., m, \mathcal{R}h^{\frac{1}{q_0}} \circ A_j \in \mathcal{A}(p_{-},q_0)$\newline\newline
 and so,
 \begin{equation*}
 \leq c \sup_{\left\Vert h\right\Vert_{\widetilde{q}'(.)} = 1} \int_{\mathbb{R}^{n}} F_{\lambda}(x) \mathcal{R}h(x) dx = c \sup_{\left\Vert h\right\Vert_{\widetilde{q}'(.)} = 1} \int_{\mathbb{R}^{n}} F_{\lambda}(x) (\mathcal{R}h^{\frac{1}{q_0}}(x))^{q_0} dx,
 \end{equation*}
 and by (\ref{DI}), since $\mathcal{R}h^{\frac{1}{q_0}} \in \mathcal{A}(p_{-},q_0)$ and $Rh \in \mathcal{A}_{1}\subset \mathcal{A}_{\infty}$,  
 \begin{equation*}
 \leq c \sup_{\left\Vert h\right\Vert_{\widetilde{q}'(.)} = 1} \sup_{\lambda >0} \int_{\mathbb{R}^n}\lambda^{q_0}\chi_{ \{ x: \sum_{i=1}^m \mathcal{M}_{\alpha}f(A^{-1}_i x)> c \lambda \} } (\mathcal{R}h^{\frac{1}{q_0}}(x))^{q_0} dx
 \end{equation*}
 Since,
 \begin{equation*}
 \left \{ x: \sum_{i=1}^m \mathcal{M}_{\alpha}f(A^{-1}_i x)> c \lambda \right \} \subseteq \bigcup_{i=1}^m \left \{ x: \mathcal{M}_{\alpha}f(A^{-1}_i x)> \frac{c \lambda}{m} \right \}
 \end{equation*}
 then,
 \begin{equation*}
  \chi_{ \{ x: \sum_{i=1}^m \mathcal{M}_{\alpha}f(A^{-1}_i x)> c \lambda \} } \leq \sum_{i=1}^m \chi_{ \{ x: \mathcal{M}_{\alpha}f(A^{-1}_i x)> \frac{c \lambda}{m} \} }.
 \end{equation*}
 so 
 \begin{equation*}
 \leq c \sup_{\left\Vert h\right\Vert_{\widetilde{q}'(.)} = 1} \sup_{\lambda >0} \sum_{i=1}^m \int_{\mathbb{R}^n}\lambda^{q_0} \chi(x)_{ \{ x: \mathcal{M}_{\alpha}f(A^{-1}_i x)> \frac{c \lambda}{m} \} } (\mathcal{R}h^{\frac{1}{q_0}}(x))^{q_0} dx
 \end{equation*} 
%If we denote by,
%\begin{equation*}
%C= \left \{ x: \mathcal{M}_{\alpha}f(A^{-1}_i x)> \frac{c \lambda}{m} \right \}
%\end{equation*}
%\begin{equation*}
%B= \left \{ y: \mathcal{M}_{\alpha}f(y)> \frac{c \lambda}{m} \right \}
%\end{equation*}
%then $x \in C$ implies that $y=A^{-1}_i x \in B$. So (\ref{E2}) is bounded, making a change of variable previous , by

\begin{equation*}
= c \sup_{\left\Vert h\right\Vert_{\widetilde{q}'(.)} = 1} \sup_{\lambda >0} \sum_{i=1}^m \int_{ \{ x: \mathcal{M}_{\alpha}f(A^{-1}_i x)> \frac{c \lambda}{m} \} }\lambda^{q_0} (\mathcal{R}h^{\frac{1}{q_0}}(x))^{q_0} dx
\end{equation*}
\begin{equation*}
= c \sup_{\left\Vert h\right\Vert_{\widetilde{q}'(.)} = 1} \sup_{\lambda >0} \sum_{i=1}^m \lambda^{q_0} \lvert det(A_i) \rvert \int_{A^{-1}_i {\left \{ x: \mathcal{M}_{\alpha}f(A^{-1}_i x)> \frac{c \lambda}{m} \right \}}} (\mathcal{R}h^{\frac{1}{q_0}}(A_i y))^{q_0} dy,
\end{equation*}
\begin{equation*}
\leq c \sup_{\left\Vert h\right\Vert_{\widetilde{q}'(.)} = 1} \sup_{\lambda >0} \sum_{i=1}^m \lambda^{q_0} \int_{\left \{ y: \mathcal{M}_{\alpha}f(y)> \frac{c \lambda}{m} \right \}} (\mathcal{R}h^{\frac{1}{q_0}}(A_i y))^{q_0} dy,
\end{equation*}
%Since ,
\begin{equation*}
\leq  c \sup_{\left\Vert h\right\Vert_{\widetilde{q}'(.)} = 1} \sup_{\lambda >0} \sum_{i=1}^m \left(\int_{\mathbb{R}^n} \lvert f(y) \rvert^{p_{-}} (\mathcal{R}h^{\frac{p_{-}}{q_0}}(A_i y)) dy \right)^{q_0},
\end{equation*}
\begin{equation*}
= c \sup_{\left\Vert h\right\Vert_{\widetilde{q}'(.)} = 1} \sum_{i=1}^m \left(\int_{\mathbb{R}^n} \lvert f(y) \rvert^{p_{-}} (\mathcal{R}h^{\frac{p_{-}}{q_0}}(A_i y)) dy \right)^{\frac{q_0}{p_{-}}},
\end{equation*}  
where the last inequality follows since $\mathcal{R}h^{\frac{1}{q_0}} \circ A_i \in \mathcal{A}(p_{-},q_0)$ for all $i=1...m$.\newline
Now we follow as in the proof of Theorem 6 to obtain
\begin{equation*}
\leq c \lVert f \rVert_{p(.)}^{q_0}.
\end{equation*}
 
\end{proof7}

Marta Urciuolo, FAMAF, UNIVERSIDAD NACIONAL DE CORDOBA, CIEM, CONICET,
Ciudad Universitaria, 5000 C\'{o}rdoba, Argentina.

E-mail adress: urciuolo@gmail.com

Lucas Vallejos, FAMAF, UNIVERSIDAD NACIONAL DE CORDOBA, CIEM, CONICET,
Ciudad Universitaria, 5000 C\'{o}rdoba, Argentina.

E-mail adress: lvallejos@famaf.unc.edu.ar

\end{document}